\documentclass[12pt]{article}
\usepackage[final]{epsfig}
\usepackage{graphics}
\usepackage{amsmath}
\usepackage{amsfonts}
\usepackage{latexsym}
\usepackage{amssymb}
\usepackage{graphicx}

\begin{document}
\newcommand{\eps}{{\varepsilon}}
\newcommand{\g}{{\gamma}}
\newcommand{\G}{{\Gamma}}
\newcommand{\proofend}{$\Box$\bigskip}
\newcommand{\C}{{\mathbf C}}
\newcommand{\Q}{{\mathbf Q}}
\newcommand{\R}{{\mathbf R}}
\newcommand{\Z}{{\mathbf Z}}
\newcommand{\RP}{{\mathbf {RP}}}
\newcommand{\CP}{{\mathbf {CP}}}
\newcommand{\Tr}{{\rm Tr\ }}
\def\proof{\paragraph{Proof.}}

\title {The equal tangents problem}
\author{Serge Tabachnikov\thanks{
Department of Mathematics,
Pennsylvania State University, University Park, PA 16802, USA;
e-mail: \tt{tabachni@math.psu.edu}
}
}

\maketitle

\section{The problem} \label{first}

Given a point $A$ outside of a closed strictly convex plane curve $\gamma$, there are two tangent segments from $A$ to $\gamma$, the left and the right ones, looking from $A$. 

Problem: {\it Does there exist a curve $\gamma$ such that one can walk around it so that, at all moments, the right tangent segment is smaller than the left one?} In other words, does there exist a pair of simple closed curves, $\gamma$ and $\Gamma$, the former strictly convex, the latter containing the former in its interior, such that for every point $A$ of $\Gamma$ the right tangent segment to $\gamma$ is smaller than the left one? 

Over the years, I have polled numerous colleagues, mostly as a dinner table topic. Most of them thought that the answer was negative, and quite a few tried to provide a proof, but each attempt had a flaw. I invite the reader to think about this question too before reading any further.

Up until recently, I believed that for any oval\footnote{a smooth strictly convex closed curve} $\gamma$ and every closed curve $\Gamma$ going around $\gamma$, there exists a point $A\in\Gamma$ from which the tangent segments to $\gamma$ are equal. In fact, I conjectured in \cite{Ta} that there existed at least {\it four} such points. (For example, this is so if $\gamma$ is an ellipse or its small perturbation).

A certain confirmation to this conjecture is the following,  lesser known and quite beautiful, theorem \cite{RMBC}, motivated by the  flotation theory. Let $\gamma$ be an oval and $\varphi$ be an angle between $0$ and $\pi$. Let $\Gamma_{\varphi}$ be the locus of points in the exterior of $\gamma$ from which $\gamma$ is seen under angle $\varphi$. Then $\Gamma_{\varphi}$ contains at least four points from which the tangent segments to $\gamma$ are equal.

In the limit $\varphi \to \pi$, the points of $\Gamma_{\varphi}$ from which the tangent segments to $\gamma$ are equal correspond to the curvature extrema of $\gamma$. This implies the famous 4-vertex theorem: a plane oval has at last four distinct curvature extrema. In the limit $\varphi \to 0$, the points of $\Gamma_{\varphi}$ from which the tangent segments to $\gamma$ are equal correspond to the diameters of $\gamma$, that is, its binormal chords, and every diameter contributes two points on the curve ``at infinity" $\Gamma_0$. This implies a well known fact that  every oval has at least two diameters, the maximal one and a minimax one, corresponding to its least width.

\begin{figure}[hbtp]
\centering
\includegraphics[width=2in]{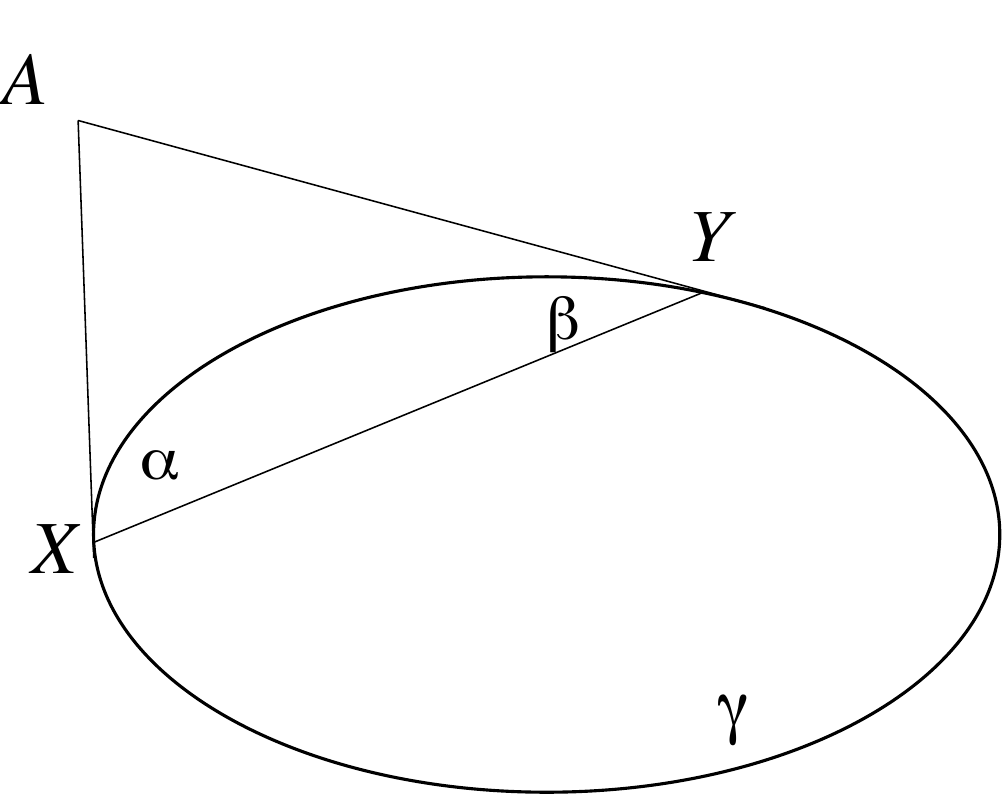}
\caption{A chord of an oval}
\label{oval}
\end{figure}

One can reformulate our problem as follows. Let $AX$ and $AY$ be the tangent segments to an oval $\gamma$, and let $\alpha$ and $\beta$ be the angles between  $XY$ and $\gamma$, see Figure \ref{oval}. The question is whether points $X$ and $Y$  can make a complete circuit along $\gamma$ so that always 
\begin{equation} \label{angles}
\beta<\alpha\ \ {\rm and}\ \  \alpha+\beta<\pi.
\end{equation}

\section{Construction} \label{second}

The goal of this note is to construct such an oval. First we construct a convex polygon $\gamma$ with a desired one-parameter family of chords, and then we approximate this polygon by an oval. Dealing with a convex polygon, we need to explain what we mean by  ``tangent lines" at its vertices. These are the support lines, that is, the lines thought  a vertex that do not intersect the interior of the polygon.

\begin{figure}[hbtp]
\centering
\includegraphics[width=5.5in]{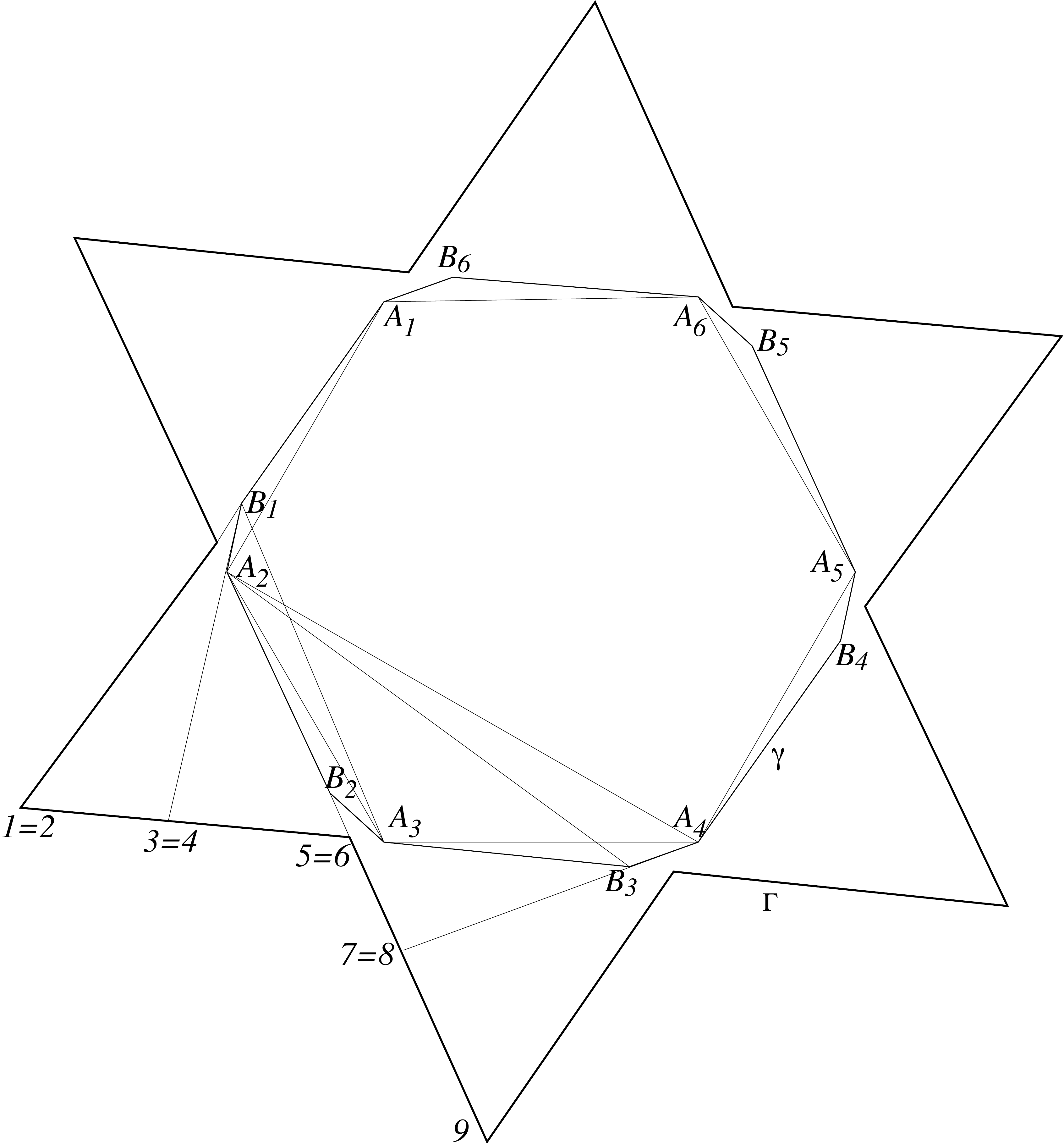}
\caption{The dodecagon $\gamma$}
\label{polygon}
\end{figure}

The polygon $\gamma$ is shown in Figure \ref{polygon}. It is a dodecagon, constructed by attaching congruent triangles, such as $A_1B_1A_2$, to the sides of a regular hexagon $A_1\dots A_6$. The point $B_i$ is sufficiently close to the  side $A_iA_{i+1}$ and closer to $A_{i+1}$ than to $A_i$. The polygon has the 6-fold rotational symmetry.

Now we describe the motion of the chord inside $\gamma$. When dealing with a smooth curve, a chord uniquely determines the tangent directions at its end-points. Since $\gamma$ is a polygon, we describe a position of a chord by a triple whose first element is the chord, the second element is a support line at the first end-point, and the third element a support line at the second end-point. 

The motion of the chord will be piecewise linear: in each step, either one end-point of the chord moves along a side of the polygon $\gamma$, or the end-points remain fixed but one of the support lines at the end-points revolves about this point. Thus, at each step, only one element of a triple changes.

Here is the whole process consisting of nine steps:

\begin{equation} \label{motion}
\begin{split}
(A_1A_3,A_1B_1,A_3B_3)\to (B_1A_3,A_1B_1,A_3B_3)\to (B_1A_3,B_1A_2,A_3B_3)\to\\
 (A_2A_3,B_1A_2,A_3B_3)\to (A_2A_3,A_2B_2,A_3B_3)\to (A_2B_3,A_2B_2,A_3B_3)\to\\ 
 (A_2B_3,A_2B_2,B_3A_4)\to (A_2A_4,A_2B_2,B_3A_4)\to (A_2A_4,A_2B_2,A_4B_4).
\end{split}
\end{equation}
After that, the process repeats using the 6-fold rotational symmetry.

The 6-pronged star $\Gamma$ around the dodecagon $\gamma$ in Figure \ref{polygon} is the locus of the intersection points of the tangent lines at the end-points of the moving chord inside $\gamma$, i.e., the locus of points labelled $A$ in Figure \ref{oval}. The points corresponding to the nine steps of the motion (\ref{motion}) are marked 1 to 9. We note that these points coincide pair-wise.

\begin{figure}[hbtp]
\centering
\includegraphics[width=2in]{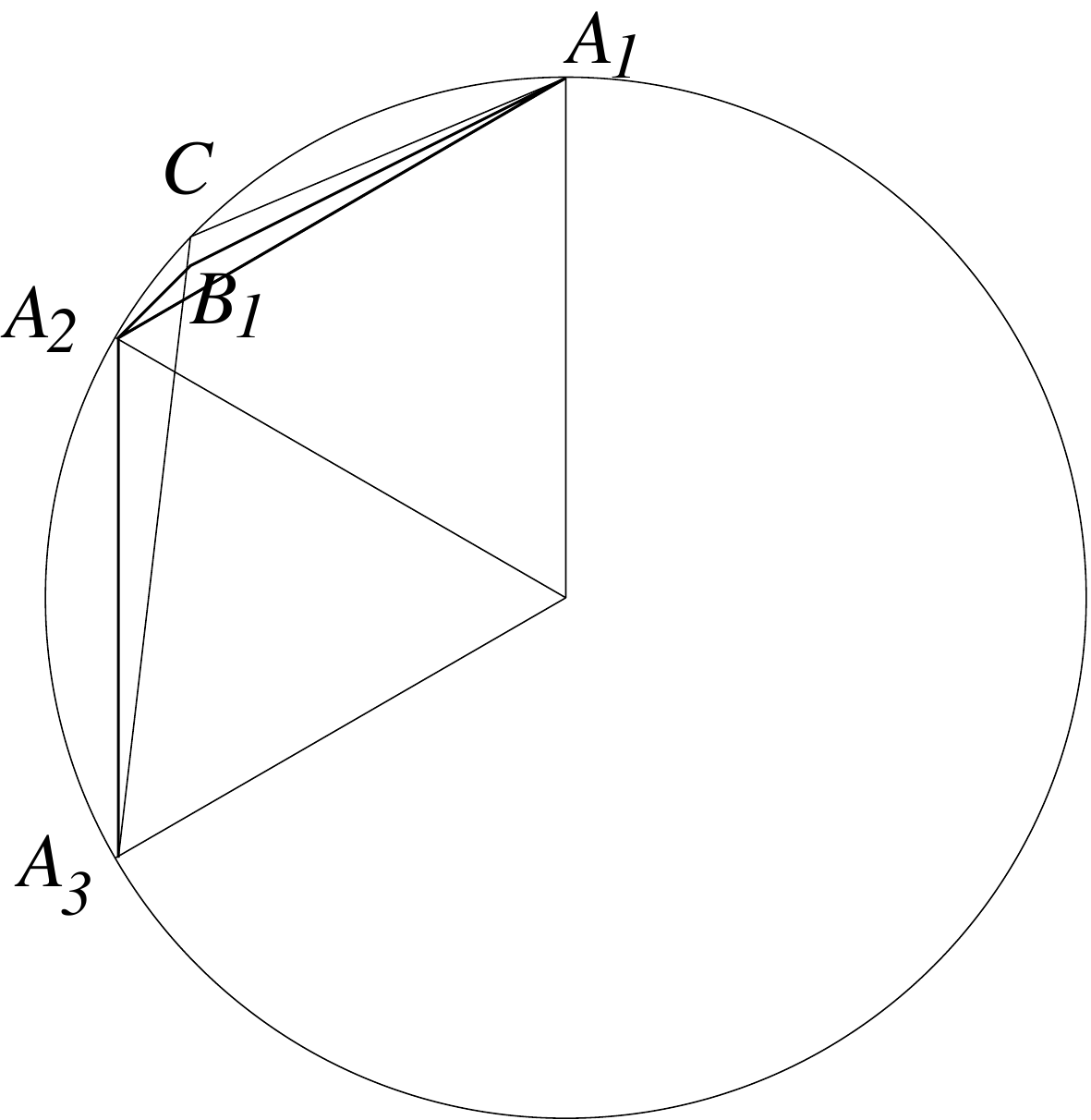}
\caption{An angle inequality}
\label{angs}
\end{figure}

We now check that the inequality $\beta < \alpha$, as in (\ref{angles}), holds during the motion (\ref{motion}). Since all the angles change monotonically, it suffices to check the inequalities at the first eight steps of (\ref{motion}). Let
$$
\angle A_2A_1B_1=\varphi,\ \angle A_1A_2B_1=\psi,\ \angle A_2A_3B_1=\theta,\ \angle B_3A_2A_4=\delta.
$$ 
We may make all these angles sufficiently small. By construction, $\varphi < \psi$. We also claim that
\begin{equation} \label{ineq}
\varphi < \theta,\ \ \varphi <\delta.
\end{equation}
Indeed, since the inscribed angles subtended by the same arc of a circle are equal, one has
$$
\theta = \angle A_2A_3C = \angle A_2A_1C > \angle A_2A_1B_1 = \varphi,
$$
see Figure \ref{angs}, and the other inequality in (\ref{ineq})  is similar.

The pairs of angles $(\beta,\alpha)$ for the first eight positions in (\ref{motion}) are easily found by elementary geometry:
\begin{equation*} 
\begin{split}
(30^{\circ}+\varphi, 90^{\circ}-\varphi),\ (60^{\circ}+\varphi-\theta, 60^{\circ}-\varphi+\theta),\ (60^{\circ}-\theta-\psi,60^{\circ}+\theta-\varphi),\\
 (60^{\circ}-\psi, 60^{\circ}-\varphi),\ (\varphi, 60^{\circ}-\varphi),\ (30^{\circ}-\delta+\varphi, 30+\delta-\varphi),\\
 (30^{\circ}-\delta+\varphi, 30^{\circ}+\psi+\delta),\ (30^{\circ}+\varphi, 30^{\circ}+\psi).
\end{split}
\end{equation*}
For each pair, the inequality $\beta<\alpha$ holds, due to the inequalities (\ref{ineq}) and $\varphi < \psi$, and the fact that  $\varphi$ and $\psi$ are small enough.

Next, we approximate the dodecagon $\gamma$ by a smooth strictly convex curve, say, $\gamma_1$. If the approximation is fine enough, 
we obtain a one-parameter family of chords of $\gamma_1$, such as depicted in Figure \ref{oval}, for which the inequality $\beta < \alpha$ still holds.
Thus $\gamma_1$ is a desired example. 

To be concrete, one may construct a $C^1$-smooth example by replacing each vertex of the dodecagon $\gamma$ by an arc of a circle of a very small radius, and every side of $\gamma$ by an arc of a circle of a very large radius. We assume that the small and the large circles share tangent directions at their common points. The resulting piecewise circular curve $\gamma_1$ is an example; see \cite{BG} concerning the fascinating geometry of piecewise circular curves.  
We note that when  $\gamma$ is approximated by, say, a piecewise circular curve, the polygon $\Gamma$ also changes slightly. In particular, the pairs of coinciding points in Figure \ref{polygon}, such as $1$ and $2$, separate and become pairs of distinct close points.

\section{Comments} \label{third}

1. Two tangent segments to an oval $\gamma$ are equal if and only if there is a circle touching $\gamma$ at  these tangency points. The locus of centers of such bi-tangent circles is called the symmetry set of $\gamma$ \cite{BGG}. Thus  our problem is closely related to  the geometry of symmetry sets. Symmetry sets has attracted much interest, in particular, due to their applications to image recognition and computer vision.
\medskip

\begin{figure}[hbtp]
\centering
\includegraphics[width=2.5in]{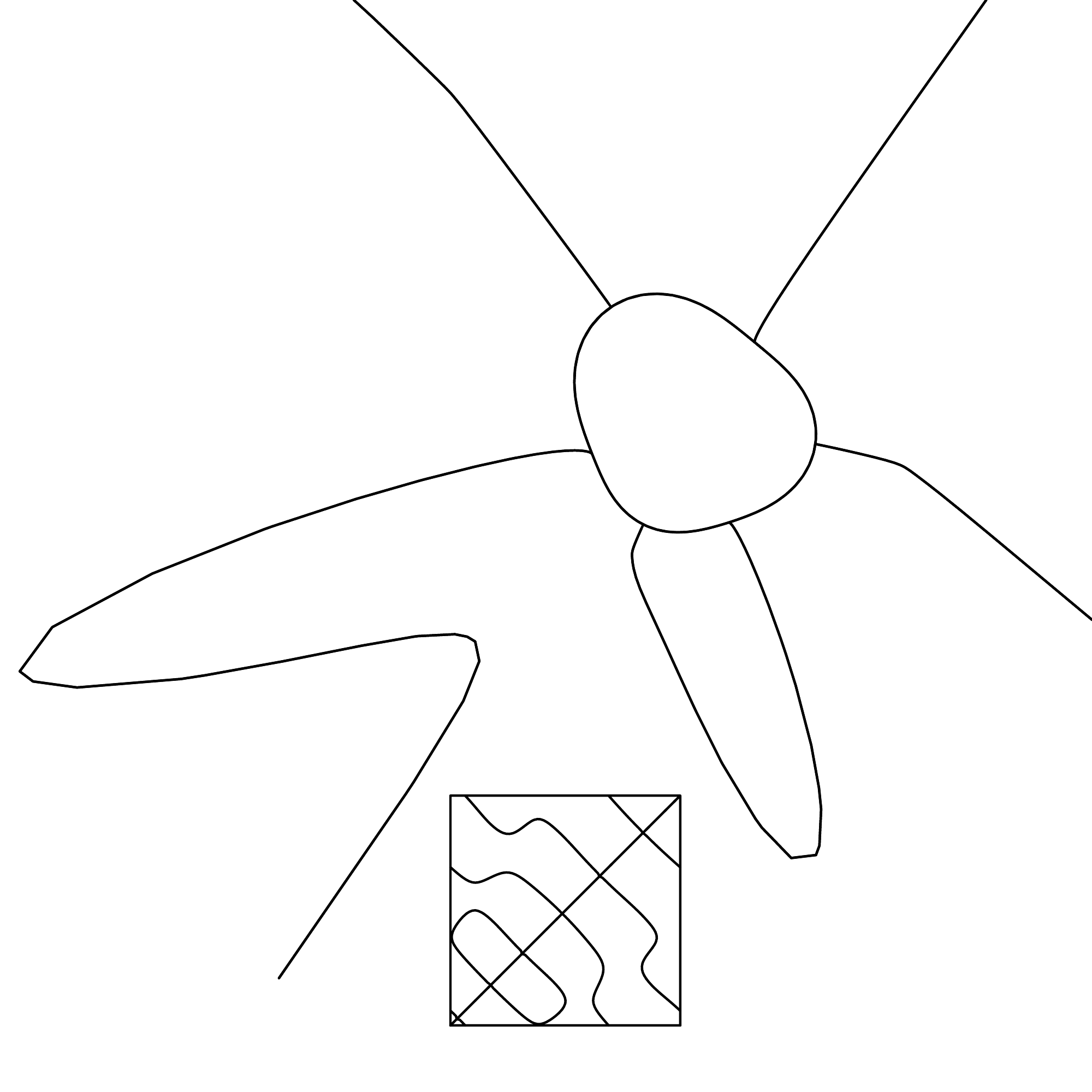}
\caption{The equitangent locus  and the respective curve on the torus}
\label{equitang}
\end{figure}

2. Given an oval $\gamma$, let $\Delta$ be the locus of points in the exterior of $\gamma$ from which the tangent segments to $\gamma$ are equal. 
Call $\Delta$ the equitangent locus of $\gamma$. 
Generically, $\Delta$ is a curve. It may have component of the following types: starting and ending on $\gamma$ (at its vertices); starting and ending at infinity; closed components; and  components starting on $\gamma$ and ending at infinity. It is this last type of components that cannot be avoided by loops going around $\gamma$. The conjecture in \cite{Ta} was that there existed at least four such components for every plane oval. See the computer-generated Figure \ref{equitang} (courtesy to P. Giblin).

3. If $\gamma$ is a convex polygon, a definition of a tangent segment to $\gamma$ from an exterior point, say, $A$, is needed. If $A$ does not lie on the extension of a side, this is a support segment to the respective vertex of the polygon. If $A$ belongs to an extension of a side, say, $XY$,  then every segment $AZ$, with point $Z$ on the side $XY$, counts as a tangent segment. Using this definition, one can define the equitangent locus $\Delta$ for a convex polygon $\gamma$; this is also a polygonal curve.
Note that, for a polygon, the relation between exterior points and chords, as depicted in Figure \ref{oval}, is not one-to-one anymore.

\begin{figure}[hbtp]
\centering
\includegraphics[width=3in]{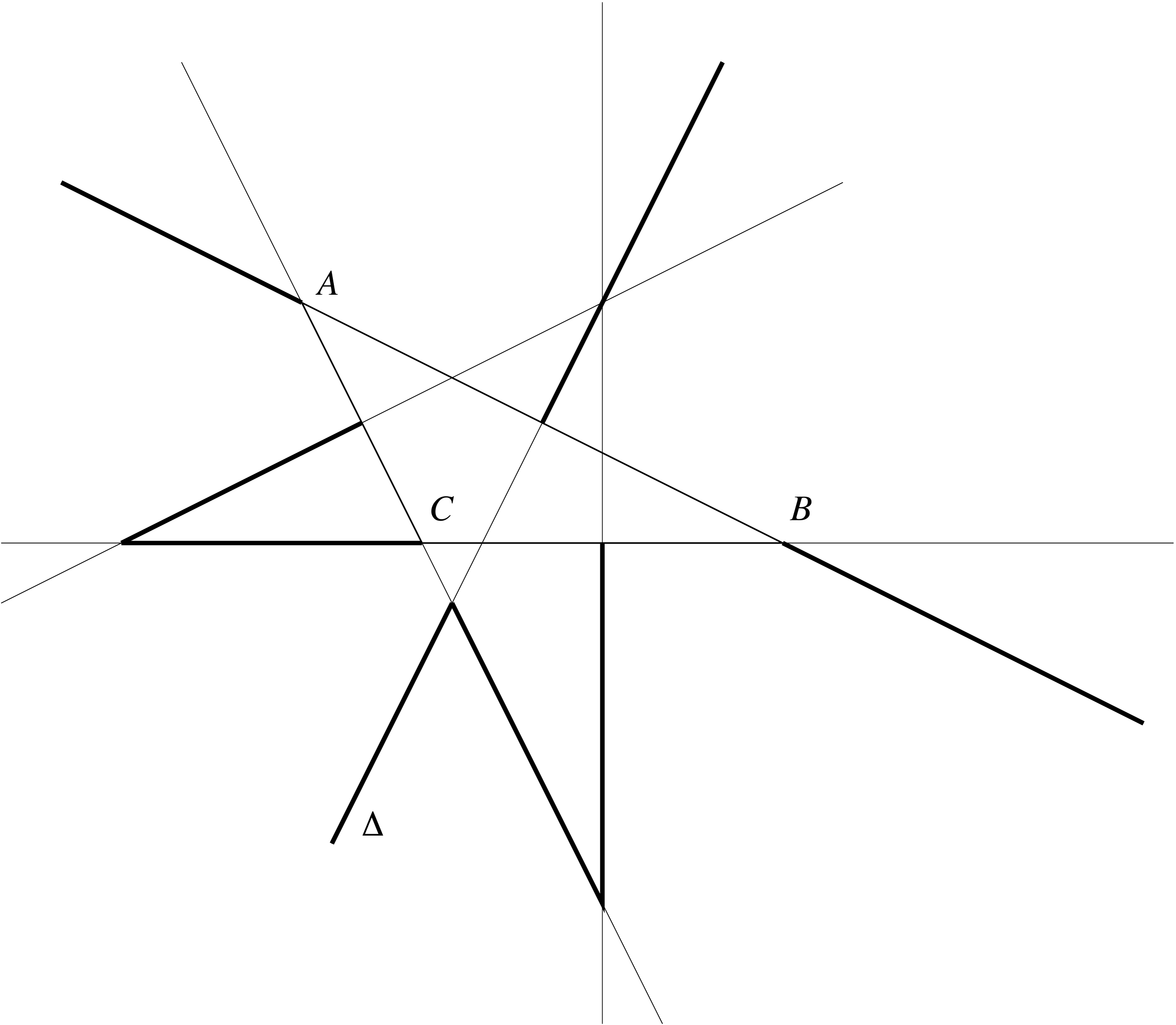}
\caption{The equitangent locus $\Delta$ for  triangle $ABC$}
\label{tri}
\end{figure}

See Figure \ref{tri} for  $\Delta$ when $\gamma$ is an obtuse  triangle. $\Delta$ is made of segments of the sides and median perpendiculars of the triangle; it has four components that start on $\gamma$ and go to infinity, and one component that starts and ends on the triangle.

4. Generically, the set of pairs of tangency points $(X,Y)$ of equal tangent segments $AX$ and $AY$ to an oval $\gamma$ is a closed curve  ${\cal C}$ on the torus $\gamma\times\gamma$. This curve is symmetric with respect to the diagonal, and the diagonal of the torus is its component. The components of ${\cal C}$ that are isotopic to the anti-diagonal are called {\it essential loops} in \cite{KO}. 

\begin{figure}[hbtp]
\centering
\includegraphics[width=2.5in]{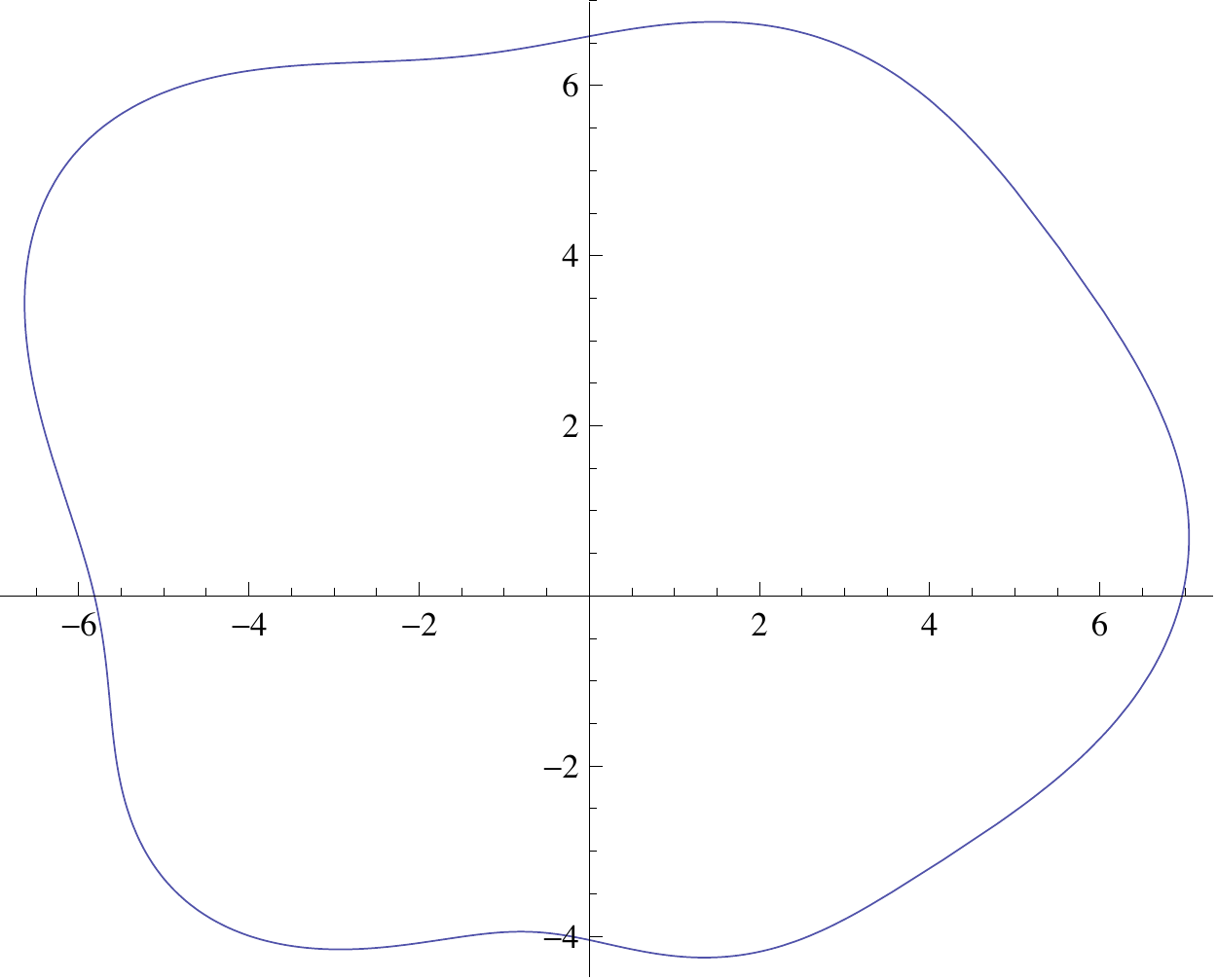}
\caption{The example of Giblin and Zakalyukin}
\label{example}
\end{figure}

A path around the oval $\gamma$ is represented by a curve on the torus that is isotopic to the diagonal. Such a curve intersects each essential loop at least twice. The example constructed in this note is free from essential loops. In 2004, P. Giblin and V. Zakalyukin constructed an example of a  non-convex plane curve free from essential loops; this example was  adapted from an earlier example of Zakalyukin \cite{Za} devised for a different purpose. See Figure \ref{example}.

\medskip

{\bf Acknowledgments}. I am grateful to all the mathematicians with whom I discussed the equal tangents problem over the years; they are too numerous to mention here by name. Special thanks go to P. Giblin and R. Schwartz for their interest and helpful discussions.


\begin{thebibliography}{99}

\bibitem{BG} T. Banchoff, P. Giblin. {\it On the geometry of piecewise circular curves.} Amer. Math. Monthly 101 (1994), 403--416.

\bibitem{BGG} J. Bruce, P. Giblin, C. Gibson. {\it Symmetry sets.} Proc. Roy. Soc. Edinburgh Sect. A 101 (1985),  163--186.

\bibitem{KO} A. Kuijper,  O. Olsen.  {\it Essential loops and their relevance for skeletons and symmetry sets}. 
Deep structure, singularities, and computer vision.  Lect, Notes  Computer Sciences, 2005,  24--35.

\bibitem{RMBC} E. ${\rm Rapha\ddot el}$, J.-M. di Meglio, M. Berger, E. Calabi. {\it Convex particles at interfaces}. 
J. Phys. I France 2 (1992) 571--579.

\bibitem{Ta} S. Tabachnikov. {\it Around four vertices.} Russian Math. Surveys 45 (1990),  229--230.

\bibitem{Za} V. Zakalyukin. {\it Maxwell stratum of Lagrangian collapse.} Proc. Steklov Inst. Math. 1998, no. 2 (221), 18--201.

\end{thebibliography}
\end{document}